\documentclass[10pt, a4paper, reqno]{amsart}
\usepackage{amscd}
\usepackage{dsfont}

\usepackage{amsmath}
\usepackage{amstext}
\usepackage{amsthm}
\usepackage{amssymb}
\usepackage{amsopn}
\usepackage{amssymb}
\usepackage{amsxtra}
\usepackage{amsfonts}
\usepackage{fancyhdr}
\usepackage{paralist, epsfig}
\usepackage[mathcal]{euscript}
\usepackage[english]{babel}
\usepackage[latin1]{inputenc}

\makeatletter
\g@addto@macro\normalsize{%
  \setlength\abovedisplayskip{10pt}
  \setlength\belowdisplayskip{10pt}
  \setlength\abovedisplayshortskip{5pt}
  \setlength\belowdisplayshortskip{8pt}
}

\allowdisplaybreaks

\newtheoremstyle{normal}
{5pt}
{5pt}
{\normalfont}
{}
{\bfseries}
{}
{0.4em}
{\bfseries{\thmname{#1}\thmnumber{ #2}.\thmnote{ \hspace{0.5em}(#3)\newline}}}

\newtheoremstyle{kursiv}
{5pt}
{5pt}
{\itshape}
{}
{\bfseries}
{}
{0.4em}
{\bfseries{\thmname{#1}\thmnumber{ #2}.\thmnote{ \hspace{0.5em}(#3)\newline}}}

\theoremstyle{kursiv}
\newtheorem{thm}{Theorem}

\theoremstyle{normal}
\newtheorem{ex}[thm]{Example}
\newtheorem{rem}[thm]{Remark}
\newtheorem*{thmA}{Desch-Schappacher Theorem}

\renewcommand{\epsilon}{\varepsilon}

\newcommand{\T}{(T(t))_{t\geqslant0}}

\newcommand{\id}{\operatorname{id}\nolimits}

\newcommand{\cs}{\operatorname{cs}\nolimits}

\renewcommand{\bar}[1]{\overline{#1}}

\usepackage{setspace}
\usepackage{color}

\definecolor{grey}{gray}{.3}

\begin{document}

\title{Desch-Schappacher perturbation of one-parameter semigroups on locally convex spaces}

\author{Birgit Jacob\hspace{1pt}\MakeLowercase{$^{\text{a,1}}$}, Sven-Ake Wegner\hspace{0.5pt}\MakeLowercase{$^{\text{b}}$} and Jens Wintermayr\hspace{0.5pt}\MakeLowercase{$^{\text{b}}$}}

\renewcommand{\thefootnote}{}
\hspace{-1000pt}\footnote{\hspace{5.5pt}2010 \emph{Mathematics Subject Classification}: Primary 47D06; Secondary 46A04, 34G10.}
\hspace{-1000pt}\footnote{\hspace{5.5pt}\emph{Key words and phrases}: Strongly continuous semigroup, Desch-Schappacher perturbation, locally convex space.\vspace{1.6pt}}

\hspace{-1000pt}\footnote{\hspace{-6.8pt}$^{\text{a},1}$\,Corresponding author: Bergische Universit\"at Wuppertal, FB C -- Mathematik, Gau\ss{}stra\ss{}e 20, 42119 Wuppertal,\linebreak\phantom{x}\hspace{12.5pt}Germany, Phone:\hspace{1.2pt}\hspace{1.2pt}+49\hspace{1.2pt}(0)\hspace{1.2pt}202\hspace{1.2pt}/\hspace{1.2pt}439\hspace{1.2pt}-\hspace{1.2pt}2527, Fax:\hspace{1.2pt}\hspace{1.2pt}+49\hspace{1.2pt}(0)\hspace{1.2pt}202\hspace{1.2pt}/\hspace{1.2pt}439\hspace{1.2pt}-\hspace{1.2pt}3724, E-Mail: bjacob@uni-wuppertal.de.\vspace{1.6pt}}

\hspace{-1000pt}\footnote{$^{\text{b}}$\,Bergische Universit\"at Wuppertal, FB C -- Mathematik, Gau\ss{}stra\ss{}e 20, 42119 Wuppertal, Germany, Phone:\hspace{1.2pt}\hspace{1.2pt}+49\hspace{1.2pt}(0)\linebreak\phantom{x}\hspace{12.5pt}202\hspace{1.2pt}/\hspace{1.2pt}439\hspace{1.2pt}-\hspace{1.2pt}2531 resp.~-2671, E-Mail: wegner@math.uni-wuppertal.de resp.~wintermayr@uni-wuppertal.de.\vspace{1.6pt}}

\begin{abstract} We prove a Desch-Schappacher type perturbation theorem for strongly continuous and locally equicontinuous one-parameter semigroups which are defined on a sequentially complete locally convex space.
\end{abstract}

\maketitle

\section{Introduction}\label{SEC-1}

Perturbation theory is an important and intensively studied topic in the context of strongly continuous operator semigroups. In the literature there exists a vast amount of results covering various types of perturbations on different classes of Banach and Hilbert spaces. For a sample we refer to the textbooks of Engel, Nagel \cite{EN}, Davies \cite{Davies} or Pazy \cite{Pazy} and the references therein. In addition, there exist important applications. Again we give only a small sample and refer to Tucsnak, Weiss \cite{TucsnakWeiss} for applications in systems theory (e.g.~closed loop systems or output feedback systems) and to Farago, Havasiy \cite{OS} for applications in numerical analysis (e.g.~operator splitting).

\medskip

Although there exists such a huge amount of results, there is no unifying theory since every type of perturbation problem requires its own methods and techniques. Below we state as an example the classical perturbation theorem of Desch, Schappacher \cite{DS}, see Engel, Nagel \cite[Corollary III.3.3]{EN}.

\medskip

\begin{thmA} Let $A\colon D(A)\rightarrow X$ be the generator of the $C_0$-semigroup $\T$ on a Banach space $X$. Let $B\colon X\rightarrow X_{-1}$ be linear and continuous. Assume that there exists $t_0>0$ such that
\begin{compactitem}
\item[(i)] $\forall\:f\in C([0,t_0],X)\colon \displaystyle\int_0^{t_0}T_{-1}(t_0-t)Bf(t)dt\in X$,
\item[(ii)] $\exists\:K\in(0,1)\;\forall\:f\in C([0,t_0],X)\colon \displaystyle\Bigl\|\int_0^{t_0}T_{-1}(t_0-t)Bf(t)dt\Bigr\|\leqslant K\textstyle\sup_{t\in[0,t_0]}|f(t)|$
\end{compactitem}\vspace{5pt}
holds. Then the perturbed operator $(A_{-1}+B)|_{X}$ again generates a $C_0$-semigroup on $X$.
\end{thmA}

\medskip

Above, the space $X_{-1}$ denotes the extrapolation space \cite[Section II.5]{EN} associated with the strongly semigroup $\T$ and $(T_{-1}(t))_{t\geqslant0}$ denotes the extension of $\T$ to $X_{-1}$. However, $X_{-1}$ can also be replaced by some more general superspace $\bar{X}\supseteq X$ as it was done by Desch and Schappacher in their original version \cite[Theorem on p.~330]{DS}.

\medskip

If we consider instead of the operator $B\colon X\rightarrow X_{-1}$ a perturbing operator $B\colon X_1\rightarrow X$, where $X_1$ denotes the interpolation space associated with the initial semigroup $\T$, then we end up with another classical type of perturbation. We refer to \cite[Section III.3]{EN} for corresponding results on so-called Miyadera-Voigt perturbations \cite{Miy66, Voi77}. In addition to the well-studied case of Banach spaces, the latter perturbations have also been considered for strongly continuous and locally equicontinuous semigroups on locally convex spaces: in 1974, Dembart \cite[Section 4]{D} proved a locally convex version of the Miyadera-Voigt perturbation theorem, see \cite[Corollary III.3.16]{EN}. In this article we close the gap and prove also a Desch-Schappacher perturbation theorem for strongly continuous and locally equicontinuous semigroups on locally convex spaces.

\medskip

For the theory of locally convex spaces we refer to Meise, Vogt \cite{MV} or Jarchow \cite{Jarchow}. For the theory of semigroups on Banach spaces we refer to Engel, Nagel \cite{EN}. Equicontinuous and exponentially equicontinuous semigroups were treated by Miyadera \cite{M}, Yosida \cite{Yosida}, Choe \cite{Choe} and recently by Albanese, Bonet, Ricker \cite{ABR10, ABR12, ABR13}. Locally exponentially semigroups were studied by K{\=o}mura \cite{Komura}, Dembart \cite{D} and recently by Albanese, K\"uhnemund \cite{AK}.

\medskip

\section{Result}\label{SEC-2}

For the whole article let $X$ be a locally convex space and denote by $cs(X)$ the set of all continuous seminorms on $X$. By a $C_0$-semigroup on $X$ we understand a family $\T$ of linear and continuous operators $T(t)\colon X\rightarrow X$ for $t\geqslant0$ such that $T(0)=\id_X$, $T(t+s)=T(t)T(s)$ for $t$, $s\geqslant0$ and $\lim_{t\rightarrow t_0}T(t)x=T(t_0)x$ for $x\in X$ and $t_0\geqslant0$ holds. The $C_0$-semigroup $\T$ is said to be locally equicontinuous if for some or, equivalently, for all $t_0>0$ the set $\{T(t)\:;\:t\in[0,t_0]\}$ is an equicontinuous subset of the space of linear and continuous maps from $X$ to itself, i.e.,
\begin{equation*}
\forall\:p\in\Gamma\;\exists\:q\in\Gamma,M\geqslant0\;\forall\:x\in X,\,t\in[0,t_0]\colon p(T(t)x)\leqslant M q(x) 
\end{equation*}
holds for some or, equivalently, for every fundamental system $\Gamma\subseteq\cs(X)$. If $\T$ is a locally equicontinuous $C_0$-semigroup, its generator is the linear operator $A\colon D(A)\rightarrow X$ given by
$$
Ax=\lim_{t\searrow0}\frac{T(t)x-x}{t} \;\text{ for }\; x\in D(A)=\bigg\{x\in X\:;\:\lim_{t\searrow0}\frac{T(t)x-x}{t}\text{ exists}\bigg\}.
$$
We refer to K\={o}mura \cite[Section 1]{Komura} and Albanese, K\"uhnemund \cite[Section 2]{AK} for the basic properties of $A\colon D(A)\rightarrow X$ and to Dembart \cite[Section 3]{D} for characterizations of those operators $A$ which appear as the generators of locally equicontinuous $C_0$-semigroups; one of his results we restate later in Theorem \ref{DEM}. Here is the statement of our main result.

\medskip

\begin{thm}\label{THM} Let $A\colon D(A)\rightarrow X$ be the generator of a locally equicontinuous $C_0$-semigroup $\T$ on a sequentially complete locally convex space $X$. Let $\bar{X}$ be a sequentially complete locally convex space such that\vspace{4pt}
\begin{compactitem}
\item[(a)] $X\subseteq\bar{X}$ is dense and the inclusion map is continuous,\vspace{4pt}
\item[(b)] $\bar{A}$ generates a locally equicontinuous $C_0$-semigroup $(\bar{T}(t))_{t\geqslant0}$ on $\bar{X}$ such that $\bar{T}(t)|_X=T(t)$ holds for all $t\geqslant0$ and with domain $D(\bar{A})=X$.\vspace{4pt}
\end{compactitem}
Let $B\colon X\rightarrow\bar{X}$ be a linear and continuous operator and $t_0>0$ be a number such that\vspace{2pt}
\begin{compactitem}
\item[(c)] $\forall\:f\in C([0,t_0],X)\colon{\displaystyle\int_0^{t_0}}\bar{T}(t_0-t)Bf(t)ds\in X$,\vspace{2pt}
\item[(d)] $\forall\:p\in\Gamma\;\exists\:K\in(0,1)\;\forall\:f\in C([0,t_0],X)\colon p\hspace{0.5pt}\Bigl({\displaystyle\int_0^{t_0}}\bar{T}(t_0-t)Bf(t)dt\Bigr)\leqslant K{\displaystyle\sup_{t\in[0,t_0]}}p(f(t))$.\vspace{3pt}
\end{compactitem}
Then the operator $C\colon D(C)\rightarrow X$ defined by
\begin{equation*}
Cf=(\bar{A}+B)x \; \text{ for } \; x\in D(C)=\big\{x\in X\:;\:(\bar{A}+B)x\in X\big\}
\end{equation*}
generates a locally equicontinuous $C_0$-semigroup on $X$ if and only if $D(C)\subseteq X$ is dense.
\end{thm}

Before we give the proof in Section \ref{SEC-3}, we state the following comments on the above results.

\begin{rem}
\begin{compactitem}
\item[(i)] In the case of Banach spaces, the prototype of a space $\bar{X}$ with the properties assumed in Theorem \ref{THM} is the extrapolation space $X_{-1}$, see Engel, Nagel \cite[Section II.5]{EN}. The latter is the natural framework for the classical perturbation theorem of Desch, Schappacher \cite{DS}, see \cite[Section III.3.a]{EN}. For exponentially equicontinuous $C_0$-semigroups $\T$ on locally convex spaces, i.e., $C_0$-semigroups satisfying
\begin{equation*}
\exists \omega \in \mathbb R \forall\:p\in\Gamma\;\exists\:q\in\Gamma,\,M\geqslant1\;\forall\:t\geqslant0,\,x\in X\colon p(T(t)x)\leqslant M e^{\omega t} q(x)
\end{equation*}
for some or, equivalently, for every fundamental system of seminorms $\Gamma\subseteq\cs(X)$, extrapolation spaces have been constructed and studied by Wegner \cite[Section 3]{W}.

\vspace{4pt}

\item[(ii)] The strategy used in \cite[Theorem III.3.1]{EN} for the Banach space proof of the Desch-Schappacher theorem is to write down explicitly a candidate for the perturbed semigroup and then to verify that its generator equals the perturbed generator of the initial semigroup. This strategy can be adapted for the locally convex case under the additional assumption that $\T$ is exponentially equicontinuous, since then Laplace transform methods are applicable. This proof is straight forward though technical and has the disadvantage that a priori it covers only a strictly smaller class of semigroups. However, in contrast to our proof below, this alternative strategy would establish en passant the variation of parameters formula
\begin{equation*}
S(t)=T(t)+\int_0^t\bar{T}(t-s)BS(s)ds
\end{equation*}
for the semigroup $(S(t))_{t\geqslant0}$ generated by $C\colon D(A)\rightarrow X$.

\vspace{4pt}

\item[(iii)] Our proof of Theorem \ref{THM} is inspired by the techniques developed by Dembart \cite{D} who proved the generation result which we restate in Theorem \ref{DEM}. Dembart used this result in order to establish a theory on Miyadera-Voigt type perturbations for locally equicontinuous $C_0$-semigroups on sequentially complete spaces.
\end{compactitem}
\end{rem}

\section{Proof}\label{SEC-3}

For the proof of Theorem \ref{THM} we need the following notation, compare \cite[Section 2]{D}. Let $X$ be a sequentially complete locally convex space and $t_0>0$. We denote by $C([0,t_0],X)$ the corresponding vector valued space of continuous functions. On the latter we consider the topology $\tau^{\infty}$ of uniform convergence given by the fundamental system $\Gamma^{\infty}=(p^{\infty})_{p\in\Gamma}$ where
\begin{equation*}
p^{\infty}(f)=\sup_{t\in[0,t_0]}p(f(t))
\end{equation*}
for $p\in\Gamma$ and $f\in C([0,t_0],X)$ whenever $\Gamma$ is a fundamental system for the topology of $X$. In addition, we consider the topology $\tau^1$ of $L^1$-convergence given by the fundamental system 
$\Gamma^1=(p^1)_{p\in\Gamma}$ where
\begin{equation*}
p^1(f)=\int_0^1p(f(t))dt
\end{equation*}
for $p\in\Gamma$ and $f\in C([0,t_0],X)$ where again $\Gamma$ is a fundamental system for the topology of $X$. We put $\mathcal{X}=\{f\in C([0,t_0],X)\:;\:f(0)=0\}$ and denote by $\mathcal{X}^{\infty}$ resp.~$\mathcal{X}^1$ the space $\mathcal{X}$ endowed with the topology induced by $(C([0,t_0],X),\tau^{\infty})$ resp.~$(C([0,t_0],X),\tau^1)$. The spaces $( C([0,t_0],X),\tau^{\infty})$ and $\mathcal{X}^{\infty}$ are sequentially complete. We denote by $C^1([0,t_0],X)$ the space of continuously differentiable vector valued functions and use
\begin{equation*}
\mathcal{D}f=f' \; \text{ for } \; f\in D(\mathcal{D})=\{f\in\mathcal{X}\cap C^1([0,t_0],X)\:;\:f'(0)=0\}
\end{equation*}
to denote the differentiation operator $\mathcal{D}\colon D(\mathcal{D})\rightarrow\mathcal{X}$. The operator $\mathcal{D}$ is closed if we consider it on the space $\mathcal{X}^{\infty}$. Let $A\colon D(A)\rightarrow X$ be a linear operator, $D(A)\subseteq X$. Then we denote by $\mathcal{A}\colon D(\mathcal{A})\rightarrow\mathcal{X}$ the operator defined via
\begin{equation*}
(\mathcal{A}f)(t)=Af(t) \; \text{ for } \: f\in D(\mathcal{A})=\{f\in\mathcal{X}\:;\:f([0,t_0])\subseteq D(A) \text{ and } t\in[0,t_0]\mapsto Af(t) \text{ is continuous}\}.
\end{equation*}
For $C\colon D(C)\rightarrow X$ as in Theorem \ref{THM} we define $\mathcal{C}\colon D(\mathcal{C})\rightarrow\mathcal{X}$ analogously and get
\begin{equation*}
\begin{array}{cc}\vspace{5pt}
(\mathcal{C}f)(t)=\bar{A}f(t)+Bf(t) \; \text{ for } \; f\in D(\mathcal{C}),\\
D(\mathcal{C})=\{f\in\mathcal{X}\:;\:\forall\:t\in[0,t_0]\colon\bar{A}f(t)+Bf(t)\in X \text{ and } t\in[0,t_0]\mapsto \bar{A}f(t)+Bf(t) \text{ is continuous}\}.
\end{array}
\end{equation*}
In the sequel we will need most of the spaces and operators defined above for $X$ also for the space $\bar{X}$ as considered in Theorem \ref{THM}. In order to distinguish the corresponding objects we will write $\bar{\mathcal{X}}$, $\bar{\mathcal{A}}$, $\bar{\mathcal{D}}$ etc.~in the latter case. For $B\colon X\rightarrow\bar{X}$ we define
\begin{equation*}
\mathcal{B}\colon\mathcal{X}\rightarrow\bar{\mathcal{X}} \; \text{ via } \; (\mathcal{B}f)(t)=Bf(t).
\end{equation*}
Now we are prepared to state Dembart's generation result for locally equicontinuous $C_0$-semigroups, which will be the essential tool for our proof of Theorem \ref{THM}. Note that we modify the formulation of Dembart's result according to his comments in \cite{D} and with respect to our purposes later on.

\begin{thm}\label{DEM}{\rm(Dembart\;\cite[Theorem 3.1]{D})} Let $X$ be a sequentially complete locally convex space and $D(A)\subseteq X$ a linear subspace. A linear operator $A\colon D(A)\rightarrow X$ is the generator of a locally equicontinuous $C_0$-semigroup if and only if\vspace{3pt}
\begin{compactitem}
\item[1.] $A$ is closed and densely defined,\vspace{2pt}
\item[2.] For some or, equivalently, for every $t_0>0$ there exists a linear operator $\mathcal{R}\colon\mathcal{X}^1\rightarrow\mathcal{X}^{\infty}$, called the generalized resolvent and being continuous by (iv) below, such that\vspace{4pt}
\begin{compactitem}
\item[(i)] $\mathcal{R}(\mathcal{D}-\mathcal{A})f=f$ holds for all $f\in D(\mathcal{D})\cap D(\mathcal{A})$,\vspace{2pt}
\item[(ii)] $\mathcal{R}f\in D(\mathcal{D})$ and $\mathcal{D}\mathcal{R}f=\mathcal{R}\mathcal{D}f$ holds for all $f\in D(\mathcal{D})$,\vspace{2pt}
\item[(iii)] $\mathcal{R}f\in D(\mathcal{A})$ and $\mathcal{A}\mathcal{R}f=\mathcal{R}\mathcal{A}f$ holds for all $f\in D(\mathcal{A})$,\vspace{2pt}
\item[(iv)] for each $p\in\Gamma$ there exist $q\in\Gamma$ and $M\geqslant0$ such that $p^{\infty}(\mathcal{R}f)\leqslant M q^1(f)$ holds for all $f\in\mathcal{X}$.\vspace{3pt}
\end{compactitem}
\end{compactitem}
The operator $\mathcal{R}$ is uniquely determined by the conditions (i)--(iv) and independent of $t_0>0$; see \cite[p.~132]{D} for details. If $A\colon D(A)\rightarrow X$ generates the locally equicontinuous $C_0$-semigroup $\T$ and $t_0>0$ is arbitrary, then the formula
\begin{equation*}
(\mathcal{R}f)(t)=\int_0^tT(t-s)f(s)ds
\end{equation*}
holds for $f\in\mathcal{X}^1$, $t\in[0,t_0]$.\hfill\qed
\end{thm}

\medskip

\textit{Proof of Theorem \ref{THM}.} We check the conditions in Theorem \ref{DEM} for the operator $C\colon D(C)\rightarrow X$. To prevent confusion we denote the generalized resolvent for $C$ by $\mathcal{R}_C$, whereas the generalized resolvents for $A$ and $\bar{A}$ will be denoted by $\mathcal{R}$ and $\bar{\mathcal{R}}$, respectively. To start with, we select $t_0>0$ and $B\colon X\rightarrow\bar{X}$ as in Theorem \ref{THM}.
\smallskip
\\1.~The operator $C$ is densely defined by our assumptions. Let $(x_{\alpha})_{\alpha}\subseteq D(C)$ be a net with $x_{\alpha}\rightarrow x\in X$ and $Cx_{\alpha}=\bar{A}x_{\alpha}+Bx_{\alpha}\rightarrow y\in X$. Since $B\colon X\rightarrow\bar{X}$ is continuous, we have $Bx_{\alpha}\rightarrow Bx$ and hence, $\bar{A}x_{\alpha}\rightarrow y-Bx$. By the closedness of $\bar{A}$ it follows that $\bar{A}x=y-Bx$ and hence, $\bar{A}x+Bx=y\in X$. Therefore, $x\in D(C)$. Consequently, $C\colon D(C)\rightarrow X$ is closed.
\smallskip
\\2.~Since $A\colon D(A)\rightarrow X$ is a generator, we can select $\mathcal{R}$ as in Theorem \ref{DEM}. In addition, $\bar{A}\colon X\rightarrow \bar{X}$ is also a generator, thus we also find $\bar{\mathcal{R}}\colon\bar{\mathcal{X}}^1\rightarrow\bar{\mathcal{X}}^{\infty}$ with the corresponding properties. According to the last part of Theorem \ref{DEM} we have
\begin{equation}\label{EQ-1}
(\bar{\mathcal{R}}f)(t) = \int_0^t\bar{T}(t-s)f(s)ds = \int_0^tT(t-s)f(s)ds = (\mathcal{R}f)(t)
\end{equation}
for $f\in\mathcal{X}$ and $t\in[0,t_0]$. We define
\begin{equation}\label{EQ-2}
\mathcal{R}_C\colon\mathcal{X}^1\rightarrow\mathcal{X}^{\infty} \; \text{ via } \; \mathcal{R}_Cf=\sum_{n=0}^{\infty}(\bar{\mathcal{R}}\mathcal{B})^n\mathcal{R}f
\end{equation}
and check that the latter is well-defined. To start with, we show that $\bar{\mathcal{R}}\mathcal{B}\colon\mathcal{X}\rightarrow\mathcal{X}$ is well-defined\label{page-label1}. Let $f\in\mathcal{X}$. We claim that $(\bar{\mathcal{R}}\mathcal{B}f)(t)\in X$ holds for every $t\in(0,t_0]$ and that $[0,t_0]\rightarrow X$, $t\mapsto (\bar{\mathcal{R}}\mathcal{B}f)(t)$ is continuous. We fix $t\in(0,t_0]$ and define
\begin{equation}\label{EQ-3}
f_t\colon[0,t_0]\rightarrow X,\;f_t(s)=\begin{cases}
        \;\hspace{22pt}0 &\text{for } \; 0\leqslant s\leqslant t_0-t, \\
        \;f(s+t-t_0) &\text{for } \; t_0-t\leqslant s\leqslant t_0.
       \end{cases}
\end{equation}
The function $f_t$ is continuous and since $f_t(0)=0$ holds, it belongs to $\mathcal{X}$. We compute
\begin{equation}\label{EQ-4}
(\bar{\mathcal{R}}\mathcal{B}f)(t) = \int_0^t\bar{T}(t-s)Bf(s)ds =  \int_{t_0-t}^{t_0}\bar{T}(t_0-s)Bf(s+t-t_0)ds = \int_0^{t_0}\bar{T}(t_0-s)Bf_t(s)ds
\end{equation}
which belongs to $X$ in view of Theorem \ref{THM}(c). It remains to check the continuity. Let $t$, $r\in[0,t_0]$ and $p\in\Gamma$. Then
\begin{equation*}
p\big((\bar{\mathcal{R}}\mathcal{B}f)(t)-(\bar{\mathcal{R}}\mathcal{B}f)(r)\big)=p\Bigl(\int_0^{t_0}\bar{T}(t_0-s)B\big(f_t(s)-f_r(s)\big)ds\Bigr)\leqslant K p^{\infty}(f_t-f_r)
\end{equation*}
holds where the equality follows from \eqref{EQ-4} and the estimate follows from Theorem \ref{THM}(d). In view of the definition of $f_t$ and $f_r$ the continuity follows by taking into account that the right shift on $\mathcal{X}^{\infty}$ is strongly continuous. This establishes our claim and $\bar{\mathcal{R}}\mathcal{B}\colon\mathcal{X}\rightarrow\mathcal{X}$ is well-defined. We have to show that the series in \eqref{EQ-2} converges. For this we need an appropriate estimate for $p^{\infty}((\bar{\mathcal{R}}\mathcal{B})^n\mathcal{R}f)$. We claim
\begin{equation}\label{EQ-5}
\forall\:p\in\Gamma\;\exists\:q\in\Gamma,\,M\geqslant0,\,K\in(0,1)\;\forall\:g\in\mathcal{X},\,n\in\mathbb{N}\colon p^{\infty}\big((\bar{\mathcal{R}}\mathcal{B})^n\mathcal{R}g\big)\leqslant K^n\,M\,q^1(g).
\end{equation}
Firstly, we prove
\begin{equation}\label{EQ-6}
\forall\:p\in\Gamma\;\exists\:K\in(0,1)\;\forall\:n\in\mathbb{N},\,f\in\mathcal{X}\colon p^{\infty}\big((\bar{\mathcal{R}}\mathcal{B})^n f\big)\leqslant K^n\,p^{\infty}(f).
\end{equation}
Let $p\in\Gamma$ be given. According to Theorem \ref{THM}(d) we select $K\in(0,1)$. Now we show by induction that
\begin{equation*}
[n]\;\;\;\;\;\forall\:f\in\mathcal{X}\colon p^{\infty}\big((\bar{\mathcal{R}}\mathcal{B})^nf\big)\leqslant K^n\,p^{\infty}(f)
\end{equation*}
is true for every $n\in\mathbb{N}$. For $n=0$ the statement is trivial. For $n\geqslant0$ we show $[n]\Rightarrow[n+1]$. Let $f\in\mathcal{X}$ be given. For $t\geqslant0$ let $f_t\in\mathcal{X}$ be defined by the formula in \eqref{EQ-3}. Using \eqref{EQ-4} and the estimate in Theorem \ref{THM}(d) we get
\begin{equation*}
p\big((\bar{\mathcal{R}}\mathcal{B}f)(t)\big)=p\Bigl(\int_0^{t_0}\bar{T}(t_0-s)Bf_t(s)ds\Bigr) \leqslant K p^{\infty}(f)
\end{equation*}
for every $t\geqslant0$ and taking the induction hypothesis into account we get
\begin{equation*}
p^{\infty}\big((\bar{\mathcal{R}}\mathcal{B})^{n+1}f\big)=p^{\infty}\big((\bar{\mathcal{R}}\mathcal{B})^{n}(\bar{\mathcal{R}}\mathcal{B}f)\big)\leqslant K^n p^{\infty}(\bar{\mathcal{R}}\mathcal{B}f)=K^n\sup_{t\in[0,t_0]}p\big((\bar{\mathcal{R}}\mathcal{B}f)(t)\big)\leqslant K^{n+1}p^{\infty}(f)
\end{equation*}
which finishes the induction and yields the condition in \eqref{EQ-6}. Now we can start our proof of \eqref{EQ-5}. Let $p\in\Gamma$ be given. Select $q\in\Gamma$ and $M\geqslant0$ as in Theorem \ref{DEM}.2(iv). Select $K\in(0,1)$ as in \eqref{EQ-6}. Let $g\in\mathcal{X}$ and $n\in\mathbb{N}$ be given. Then the estimates in \eqref{EQ-6}, with $f=\mathcal{R}g\in\mathcal{X}$, and Theorem \ref{DEM}.2(iv) yield
\begin{equation*}
p^{\infty}\big((\bar{\mathcal{R}}\mathcal{B})^n\mathcal{R}g\big)\leqslant K^n p^{\infty}(\mathcal{R}g)\leqslant K^n\,M\,q^1(g)
\end{equation*}
which is the estimate in \eqref{EQ-5} and thus establishes the claim. Now we turn back to the series in \eqref{EQ-2}. We fix $f\in\mathcal{X}$ and abbreviate the partial sums of the series in \eqref{EQ-2} by
\begin{equation*}
S_jf=\sum_{n=0}^{j}(\bar{\mathcal{R}}\mathcal{B})^n\mathcal{R}f.
\end{equation*}
Using \eqref{EQ-5}, we get
\begin{equation*}
\forall\:p\in\Gamma\;\exists\:q\in\Gamma\,M\geqslant0,\,K\in(0,1)\;\forall\:j>i\colon p^{\infty}(S_jf-S_if)\leqslant \sum_{n=i+1}^jp^{\infty}\big((\bar{\mathcal{R}}\mathcal{B})^n\mathcal{R}f\big)\leqslant Mq^1(f)\sum_{n=i+1}^jK^n
\end{equation*}
which shows that $(S_jf)_{j\in\mathbb{N}}\subseteq\mathcal{X}^{\infty}$ is a Cauchy sequence. Since $\mathcal{X}^{\infty}$ is sequentially complete, $(S_jf)_{j\in\mathbb{N}}$ is convergent and thus $\mathcal{R}_C$ is well-defined.

\bigskip

In order to complete the proof we have to verify the conditions in Theorem \ref{DEM}.2(i)--(iv) with $\mathcal{A}$ replaced with $\mathcal{C}$ and $\mathcal{R}$ replaced with $\mathcal{R}_C$.

\medskip

(iv)~Let $p\in\Gamma$ be given. Select $q\in\Gamma$, $M\geqslant0$ and $K\in(0,1)$ as in \eqref{EQ-5}. Let $g\in\mathcal{X}$ be given. For $N\in\mathbb{N}$ we have
$$
p^{\infty}\Bigl(\sum_{n=0}^N(\bar{\mathcal{R}}\mathcal{B})^n\mathcal{R}_Ag\Bigr)\leqslant \sum_{n=0}^Np^{\infty}\big((\bar{\mathcal{R}}\mathcal{B})^n\mathcal{R}_Ag\big) \leqslant M q^1(g)\sum_{n=0}^N K^n \leqslant M\,\frac{1}{1-K}\,q^{1}(g)
$$
which proves condition (iv).

\medskip

(i)~Let $f\in D(\mathcal{D})\cap D(\mathcal{C})$. Thus, $(\mathcal{D}-\mathcal{C})f\in\mathcal{X}$ holds by the definitions of $D(\mathcal{C})$ and $D(\mathcal{D})$. We have
\begin{equation*}
\mathcal{R}_C(\mathcal{D}-\mathcal{C})f = \sum_{n=0}^{\infty}(\bar{\mathcal{R}}\mathcal{B})^n\mathcal{R}(\mathcal{D}-(\bar{\mathcal{A}}+\mathcal{B}))f=\sum_{n=0}^{\infty}(\bar{\mathcal{R}}\mathcal{B})^n\bar{\mathcal{R}}((\mathcal{D}-\bar{\mathcal{A}})f-\mathcal{B}f)
\end{equation*}
by the definition of $\mathcal{R}_C$ and since $\bar{\mathcal{R}}|_{\mathcal{X}}=\mathcal{R}$ holds in view of \eqref{EQ-1}. We claim
\begin{equation}\label{EQ-7}
\mathcal{R}_C(\mathcal{D}-\mathcal{C})f=\sum_{n=0}^{\infty}(\bar{\mathcal{R}}\mathcal{B})^n\bar{\mathcal{R}}(\mathcal{D}-\bar{\mathcal{A}})f - \sum_{n=0}^{\infty}(\bar{\mathcal{R}}\mathcal{B})^n\bar{\mathcal{R}}\mathcal{B}f
\end{equation}
where both series belong to $\mathcal{X}^{\infty}$ and are convergent in this space---we notice that we a priori only know that the series both belong to $\bar{\mathcal{X}}$. Let us start with the second series. On p.~\ref{page-label1} we proved that $\bar{\mathcal{R}}\mathcal{B}\colon\mathcal{X}\rightarrow\mathcal{X}$ is well-defined. Therefore, the second series in \eqref{EQ-7} is indeed a series in $\mathcal{X}$. We denote the corresponding partial sums by $T_jf$ and deduce
\begin{equation*}
\forall\:p\in\Gamma\;\exists\:K\in(0,1)\;\forall\:j>i\colon p^{\infty}(T_jf-T_if)\leqslant \sum_{n=i+1}^jp^{\infty}\big((\bar{\mathcal{R}}\mathcal{B})^nf\big)\leqslant p^{\infty}(f)\sum_{n=i+1}^jK^n
\end{equation*}
from \eqref{EQ-6}. This shows that $(T_j)_{j\in\mathbb{N}}$ is Cauchy and hence convergent in $\mathcal{X}^{\infty}$. Now we observe that
\begin{equation}\label{EQ-8}
\bar{\mathcal{R}}(\mathcal{D}-\bar{\mathcal{A}})f=\bar{\mathcal{R}}(\bar{\mathcal{D}}-\bar{\mathcal{A}})f=f
\end{equation}
holds. The first equality is valid since $\bar{\mathcal{D}}|_{D(\mathcal{D})}=\mathcal{D}$ and $D(\mathcal{D})\subseteq D(\bar{\mathcal{D}})$ follows from the continuity of the inclusion map $X\subseteq\bar{X}$. The second is the condition in Theorem \ref{DEM}.2(i) for the generator $\bar{A}$. It follows that the first series of \eqref{EQ-7} belongs to $\mathcal{X}^{\infty}$ and that it is convergent there. This establishes the claim and enables us to derive 
\begin{equation*}
\mathcal{R}_C(\mathcal{D}-\mathcal{C})f = \sum_{n=0}^{\infty}(\bar{\mathcal{R}}\mathcal{B})^nf - \sum_{n=1}^{\infty}(\bar{\mathcal{R}}\mathcal{B})^nf = f
\end{equation*}
from \eqref{EQ-7} and \eqref{EQ-8}, which consequently proves condition (i).

\bigskip

(ii)~For $0<h<t_0$, $t\in [0,t_0-h]$ and $f\in D(\mathcal{D})$ we compute
\begin{equation*}
\int_{0}^{t}\bar{T}(t-s)B\frac{1}{h}(f(s+h)-f(s))ds = \frac{1}{h}\big((\bar{\mathcal{R}}\mathcal{B}f)(t+h)-(\bar{\mathcal{R}}\mathcal{B}f)(t)\big)-\frac{1}{h}\int_{0}^{h}\bar{T}(t+h-s)Bf(s)ds.
\end{equation*}
Using Theorem \ref{THM}(d) and $f\in D(\mathcal{D})$ we get that the left hand side converges as $h\rightarrow 0$  in $X$ to $(\bar{\mathcal{R}}\mathcal{B}\mathcal{D}f)(t)$. Again using Theorem \ref{THM}(d), $f\in D(\mathcal{D})$, in particular the fact that $f(0)=f'(0)=0$, implies that  
\begin{equation*}
\frac{1}{h}\int_{0}^{h}\bar{T}(t+h-s)Bf(s)ds \rightarrow 0 \quad \mbox{in $X$ as $h\rightarrow 0$.}
\end{equation*}
Thus we obtain that $\bar{\mathcal{R}}\mathcal{B}f \in  D(\mathcal{D})$ and
 $\bar{\mathcal{R}}\mathcal{B}\mathcal{D}f=\mathcal{D}\bar{\mathcal{R}}\mathcal{B}f$. Furthermore for all $f\in D(\mathcal{D})$
\begin{equation*}
\mathcal{D}\sum_{n=0}^m(\bar{\mathcal{R}}\mathcal{B})^n\mathcal{R}f=\sum_{n=0}^m(\bar{\mathcal{R}}\mathcal{B})^n\mathcal{R}\mathcal{D}f
\end{equation*}
holds and as $\mathcal{D}$ is a closed operator we conclude
\begin{equation*}
\mathcal{R}_C\mathcal{D}f=\mathcal{D}\mathcal{R}_Cf
\end{equation*}
for all $f\in D(\mathcal{D})$.

\bigskip

(iii)~We start this last part by showing that
\begin{equation}\label{EQ-9}
\forall\:f\in D(\mathcal{D})\cap D(\mathcal{C})\colon \mathcal{R}_Cf\in D(\mathcal{C}) \: \text{ and } \: \mathcal{R}_C\mathcal{C}f=\mathcal{C}\mathcal{R}_Cf
\end{equation}
holds. Thus, let $f\in D(\mathcal{D})\cap D(\mathcal{C})$. We apply condition (i), which we proved already, and get that $\mathcal{R}_C\mathcal{C}f=\mathcal{R}_C\mathcal{D}f-f$ holds. By (ii), which we also already showed, we can permute $\mathcal{R}_C$ and $\mathcal{D}$ in the latter equation. By adding a zero and since $\mathcal{D}=\bar{\mathcal{D}}|_{\mathcal{X}}$ and $\mathcal{R}=\bar{\mathcal{R}}|_{\mathcal{X}}$ hold we get
\begin{equation}\label{EQ-10}
\mathcal{R}_C\mathcal{C}f = \mathcal{R}_C\mathcal{D}f-f = \mathcal{D}\mathcal{R}_Cf - \mathcal{D}\mathcal{R}f + \mathcal{D}\mathcal{R}f - f = \bar{\mathcal{D}}(\mathcal{R}_C-\mathcal{R})f+\bar{\mathcal{D}}\,\bar{\mathcal{R}}f-f.
\end{equation}
We compute
\begin{equation}\label{EQ-11}
(\mathcal{R}_C-\mathcal{R})f = \sum_{n=0}^{\infty}(\bar{\mathcal{R}}\mathcal{B})^n\mathcal{R}f-\mathcal{R}f = \sum_{n=1}^{\infty}(\bar{\mathcal{R}}\mathcal{B})^n\mathcal{R}_Af = \bar{\mathcal{R}}\mathcal{B}\sum_{n=0}^{\infty}(\bar{\mathcal{R}}\mathcal{B})^n\mathcal{R}f = \bar{\mathcal{R}}\mathcal{B}\mathcal{R}_Cf
\end{equation}
and obtain from \eqref{EQ-10} that
\begin{equation}\label{EQ-12}
\mathcal{R}_C\mathcal{C}f = \bar{\mathcal{D}}(\mathcal{R}_C-\mathcal{R})f+\bar{\mathcal{D}}\,\bar{\mathcal{R}}f-f= \bar{\mathcal{D}}\,\bar{\mathcal{R}}\mathcal{B}\mathcal{R}_Cf+\bar{\mathcal{D}}\,\bar{\mathcal{R}}f-f
\end{equation}
holds. Next we claim
\begin{equation}\label{EQ-13}
\forall\:g\in D(\bar{\mathcal{D}})\colon \bar{\mathcal{R}}g\in D(\bar{\mathcal{A}}) \; \text{ and } \; \bar{\mathcal{D}}\,\bar{\mathcal{R}}g = \bar{\mathcal{A}}\,\bar{\mathcal{R}}g+g.
\end{equation}
For $g\in D(\bar{\mathcal{D}})$ and $h>0$ we compute
\begin{equation*}
\frac{1}{h}[\bar{T}(h)-I]\big( (\bar{\mathcal{R}}g)(t)\big) = \frac{1}{h}\big( (\bar{\mathcal{R}}g)(t+h) - (\bar{\mathcal{R}}g)(t)\big) - \frac{1}{h}\int_0^h\bar{T}(s)g(t+h-s)ds
\end{equation*}
and let $h\searrow0$. By Theorem \ref{DEM}.2(ii), applied for $\bar{A}$, we get $\bar{\mathcal{R}}g\in D(\bar{\mathcal{D}})$. Therefore, the first summand on the right hand side converges to $(\bar{\mathcal{D}}\,\bar{\mathcal{R}}g)(t)$ in $\bar{X}$. The second summand converges in $\bar{X}$ to $g(t)$ since $(\bar{T}(t))_{t\geqslant0}$ is locally equicontinuous on $\bar{X}$. Consequently, the limit for $h\searrow0$ of the left hand side exists in $\bar{X}$, i.e., $(\bar{\mathcal{R}}g)(t)\in D(\bar{A})$ and $t\mapsto \bar{A}( (\bar{\mathcal{R}}g)(t) \big)=(\bar{\mathcal{D}}\,\bar{\mathcal{R}}g)(t)-g(t)$ is continuous since $\bar{\mathcal{D}}\,\bar{\mathcal{R}}g$, $g\in\bar{X}$. This yields $\bar{\mathcal{R}}g\in D(\bar{\mathcal{A}})$ and $\bar{\mathcal{A}}\,\bar{\mathcal{R}}g=\bar{\mathcal{D}}\,\bar{\mathcal{R}}g-g$ which establishes \eqref{EQ-13}. 

\smallskip

Now we apply \eqref{EQ-13} to \eqref{EQ-12}, for the second summand with $g=f\in D(\mathcal{D})\subseteq D(\bar{\mathcal{D}})$ and for the first summand with $g=\mathcal{B}\mathcal{R}_Cf$. The latter belongs to $D(\bar{\mathcal{D}})$ since $\mathcal{R}_C\in D(\mathcal{D})$ by (ii), which we showed already, and since $B\colon X\rightarrow\bar{X}$ is continuous. We obtain
\begin{equation}\label{EQ-14}
\mathcal{R}_C\mathcal{C}f
= \bar{\mathcal{D}}(\mathcal{R}_C-\mathcal{R})f+\bar{\mathcal{D}}\,\bar{\mathcal{R}}f-f
= \bar{\mathcal{D}}\,\bar{\mathcal{R}}\mathcal{B}\mathcal{R}_Cf+\bar{\mathcal{A}}\,\bar{\mathcal{R}}f
= \bar{\mathcal{A}}\,\bar{\mathcal{R}}\mathcal{B}\mathcal{R}_Cf+\mathcal{B}\mathcal{R}_Cf+\bar{\mathcal{A}}\mathcal{R}f.
\end{equation}
Finally we use \eqref{EQ-11} again to get
\begin{equation*}
\mathcal{R}_C\mathcal{C}f
=
\bar{\mathcal{A}}\,\bar{\mathcal{R}}\mathcal{B}\mathcal{R}_Cf+\mathcal{B}\mathcal{R}_Cf+\bar{\mathcal{A}}\mathcal{R}f
=
\bar{\mathcal{A}}(\bar{\mathcal{R}}\mathcal{B}\mathcal{R}_C+\mathcal{R})f+\mathcal{B}\mathcal{R}_Cf
= \bar{\mathcal{A}}\mathcal{R}_Cf+\mathcal{B}\mathcal{R}_Cf
= \mathcal{C}\mathcal{R}_Cf
\end{equation*}
from \eqref{EQ-14}, which establishes \eqref{EQ-9}.

\smallskip

For the final step of the proof we need the following density statement. We consider the right shift on the space $\mathcal{X}^{\infty}$, where it defines an equicontinuous $C_0$-semigroup with generator $\mathcal{D}\colon D(\mathcal{D})\rightarrow\mathcal{X}$. We put $D:=D(\mathcal{D})\cap D(\mathcal{C})$. We observe that $D(\mathcal{C})$ as well as $D(\mathcal{D})$ are invariant under the right shift. Consequently the same is true for $D$. By part 1 of this proof and \cite[Proposition 2.2]{D} it follows that $D\subseteq\mathcal{X}^{\infty}$ is dense. By \cite[Proposition 7]{AK} we get that $D\subseteq(D(\mathcal{D}),\tau_{D(\mathcal{D})})$ is dense where $\tau_{D(\mathcal{D})}$ denotes the graph topology given by the fundamental system of seminorms $\Gamma_{D(\mathcal{D})}=\{p^{\infty}(\cdot)+p^{\infty}(\mathcal{D}(\cdot))\:;\:p\in\Gamma\}$.

\smallskip

Now we are ready to prove the condition in Theorem \ref{DEM}.2(iii) for $\mathcal{C}$ and $\mathcal{R}_C$. Let $f\in D(\mathcal{C})$. Due to the last paragraph we can select a net $(f_{\alpha})_{\alpha}\subseteq D$ with $f_{\alpha}\rightarrow f$ with respect to $\tau_{D(\mathcal{D})}$. Since the canonical maps $(D(\mathcal{D}),\tau_{D(\mathcal{D})})\rightarrow \mathcal{X}^{\infty}\rightarrow \mathcal{X}^{1}$ are continuous, $f_{\alpha}\rightarrow f$ also holds in $\mathcal{X}^1$. By (iv), which we proved already, it follows that $\mathcal{R}_Cf_{\alpha}\rightarrow\mathcal{R}_Cf$ converges in $\mathcal{X}^{\infty}$. From \eqref{EQ-9} it follows $\mathcal{R}_Cf_{\alpha}\in D(\mathcal{C})$ and together with (ii), which we also proved already, we get
\begin{equation*}
\mathcal{C}\mathcal{R}_Cf_{\alpha}=\mathcal{R}_C\mathcal{C}f_{\alpha}=\mathcal{R}_C\mathcal{D}f_{\alpha}-f_{\alpha}
\end{equation*}
for every $\alpha$. Since the maps
\begin{equation*}
(D(\mathcal{D}),\tau_{D(\mathcal{D})})\stackrel{\mathcal{D}}{\longrightarrow}\mathcal{X}^{\infty}\stackrel{I}{\longrightarrow}\mathcal{X}^{1}\stackrel{\mathcal{R}_C}{\longrightarrow}\mathcal{X}^{\infty}
\end{equation*}
are all continuous, we conclude that $\mathcal{C}\mathcal{R}_Cf_{\alpha}\rightarrow g\in\mathcal{X}$ converges in the topology of $\mathcal{X}^{\infty}$. By \cite[Proposition 2.3]{D} the operator $\mathcal{C}$ is closed with respect to the topology of $\mathcal{X}^{\infty}$ and consequently $\mathcal{R}_Cf\in D(\mathcal{C})$ and
\begin{equation*}
\mathcal{C}\mathcal{R}_Cf_{\alpha}=\mathcal{R}_C\mathcal{C}f_{\alpha}
\end{equation*}
is valid.\hfill\qed

\section{Example}

The following example is a typical application of Desch-Schappacher perturbation results to boundary perturbations \cite[Example III.3.5]{EN} but in this case considered on the Fr\'{e}chet space $L^p_{\text{loc}}(-\infty,b]$ for $b>0$.

\begin{ex} Let $b>0$, $1<p<\infty$, $X=L^p_{\text{loc}}(-\infty,b]$ with fundamental system $\Gamma=\{p_n\:;\:n\in\mathbb{N}\}$ where $p_n(x)=(\int_{-n}^b|x(t)|^p\mathrm{d}t)^{1/p}$ for $n\in\mathbb{N}$ and $x\in X$. We define the operator $C\colon D(C)\rightarrow X$ as the derivative
\begin{equation*}
Cx=x' \; \text{ for } \; x\in D(C)=\big\{x\in W^{1,p}_{\text{loc}}(-\infty,b]\:;\: x(b)=\Phi(x)\big\}
\end{equation*}
where $\Phi\in X'$. In particular, there exists a number $K>0$ s.th. $|\Phi(x)|\leqslant K p_1(x)$ holds for all $x\in X$. Now we show that the operator $C$ generates an exponentially equicontinuous $C_0$-semigroup on $X$. Therefore we consider $C$ as a perturbation of the operator
\begin{equation*}
Ax=x' \; \text{ for } \; x\in D(A)=\big\{x\in W^{1,p}_{\text{loc}}(-\infty,b] \:;\:x(b)=0\big\}.
\end{equation*}

\medskip

First we show that $A$ is a generator of the left shift semigroup $\T$ defined by $(T(t)x)(s)=x(s+t)$ for $t+s\leqslant b$ and $(T(t)x)(s)=0$ for $s+t>b$. Clearly $\T$ satisfies the evolution property. We know that $(L^p[-n,b],p_n)$ is a Banach space for every $n\in\mathbb{N}$ and the restriction of $\T$ to $L^p[-n,b]$ is a $C_0$-semigroup with generator $A|_{L^p[-n,b]}$ and so $\T$ is strongly continuous at zero. Furthermore we have
\begin{equation*}
p_n(T(t)x)^p = \int_{-n}^{b-t}|x(t+s)|^p ds + \int_{b-t}^{b}|0|^p ds = \int_{-n+t}^{b}|x(s)|^p ds \leqslant \int_{-n}^b|x(s)|^p ds = p_n(x)^p
\end{equation*}
for all $n\in\mathbb{N}$, $t\geqslant0$ and $x\in X$. Therefore, $\T$ is even equicontinuous.

\medskip

Now we define the perturbation $B\colon X\rightarrow X_{-1}$ via $Bx = -\Phi(x)A_{-1}\mathds{1}$ where $\mathds{1}=\mathds{1}_{(-\infty,b]}$ denotes the characteristic function of $(-\infty,b]$ and $X_{-1}$ denotes the extrapolation space associated with $X$ and $\T$, see \cite[Section 3]{W}. We get
\begin{equation*}
(A_{-1}+B)x=A_{-1}x-\Phi(x)A_{-1}\mathds{1}=A_{-1}(x-\Phi(x)\mathds{1})=Cx
\end{equation*}
for every $x\in X$ and thus $(A_{-1}+B)|_X=C$ where $A_{-1}$ is the generator of the extended semigroup $(T_{-1})_{t\geqslant0}$ on $X_{-1}$.

\medskip

If we show that for the perturbation $B$ and some $t_0>0$ the conditions in Theorem \ref{THM} are fulfilled, it follows that $C$ is a generator of an exponentially equicontinuous $C_0$-semigroup.

\medskip

Let $f\in L^p([0,b],X)$. Then we have
\begin{eqnarray*}
\int_0^{b}T_{-1}(b-s)Bf(s)ds=\int_0^{b}T_{-1}(b-s)(-1)\Phi(f(s))A_{-1}\mathds{1} ds =- A_{-1}\int_0^{b}T(b-s)\Phi(f(s))\mathds{1}ds
\end{eqnarray*}
which belongs to $X$ if and only if
\begin{eqnarray*}
\int_0^{b}T(b-s)(\Phi(f(s))\mathds{1} ds\in D(A).
\end{eqnarray*}
We define the function
\begin{equation*}
g(\cdot)=\int_0^{b}T(b-s)(\Phi(f(s))\mathds{1}\mathrm{d}s=\int_0^{b}(\Phi(f(s))\mathds{1}(\cdot+b-s)ds
\end{equation*}
and in view of the definition of the characteristic function we see
\begin{equation}\label{EQ-EX}
g(t)=
\begin{cases}
\int_t^b\Phi(f(s))\mathrm{d}s & \text{for } t\in [0,b],\\
\int_0^b\Phi(f(s))\mathrm{d}s & \text{for } t\in (-\infty,0).
\end{cases}
\end{equation}
As $\Phi\in X'$ and $f\in L^p([0,b],X)$ it follows $\Phi\circ f\in L^p([0,b])$ and therefore $g\in W^{1,p}((-\infty,b])$ with $g(b)=0$. This shows $g\in D(A)$.

\medskip

Now let $f\in C([0,t_0],X)$ for some $t_0\in (0,b)$. We define $\tilde{f}\in L^p([0,b],X)$ via
\begin{equation*}
\tilde{f}(s)=
\begin{cases}
\;\hspace{22pt} 0 & \text{for } 0\leqslant s\leqslant b-t_0, \\
\;f(s+t_0-b) & \text{for } b-t_0\leqslant s\leqslant b.
\end{cases}
\end{equation*}
Then we get $\int_0^{t_0}T_{-1}(t_0-s)Bf(s)\mathrm{d}s=\int_0^{b}T_{-1}(b-s)B\tilde{f}(s)\mathrm{d}s\in X$ and we define the function
\begin{equation*}
h(\cdot)=\int_0^{t_0}T_{-1}(t_0-s)Bf(s)\mathrm{d}s=\int_0^{b}T_{-1}(b-s)B\tilde{f}(s)\mathrm{d}s.
\end{equation*}
From \eqref{EQ-EX} we see
\begin{equation*}
h(t)=-A\left\{\hspace{-3pt}
\begin{array}{ll}
\int_t^b\Phi(\tilde{f}(s))\mathrm{d}s &\text{for } t\in [0,b] \\
\int_0^b\Phi(\tilde{f}(s))\mathrm{d}s &\text{for } t\in (-\infty,0)
\end{array}
\hspace{-3pt}\right\}
=\left\{\hspace{-3pt}
\begin{array}{ll}
\Phi(f(t_0-b+t)) &\text{for } t\in [b-t_0,b] \\
\hspace{28pt}0 &\text{for } t\in (-\infty,b-t_0)
\end{array}
\hspace{-3pt}\right\}.
\end{equation*}
It follows
\begin{equation*}
p_n(h)^p=\int_{b-t_0}^b|\Phi(f(t_0-b+t))|^p dt\leqslant \int_{b-t_0}^bK^p p_1(f(t_0-b+t))^p dt\leqslant \,t_0\hspace{0.8pt}K^p\hspace{0.8pt}p_1^{\infty}(f)^p
\end{equation*}
and thus
\begin{equation*}
p_n\Bigl(\int_0^bT_{-1}(b-s)Bf(s) ds\Bigr) = p_n(h) \leqslant t_0^{1/p}\hspace{0.8pt}K\hspace{0.8pt}p_1^{\infty}(f)\leq t_0^{1/p}Kp_n^{\infty}(f).
\end{equation*}
Finally it is enough to select at the beginning $t_0>0$ such that $t_0^{1/p}K\in(0,1)$ holds.
\end{ex}

\bigskip

\normalsize

\bibliographystyle{amsplain}

\bibliography{DeschSchappacher}
\end{document}